\documentclass[12pt]{article}

\usepackage{amsmath,amssymb,amsbsy,amsfonts,amsthm,latexsym,
                   amsopn,amstext,amsxtra,euscript,amscd}

\newtheorem{theorem}{Theorem}[section]
\newtheorem{lemma}[theorem]{Lemma}

\theoremstyle{definition}

\def\fl#1{\left\lfloor#1\right\rfloor}
\def\rf#1{\left\lceil#1\right\rceil}

\def\cB{\mathcal B}

\def\cI{\mathcal I}

\def\Z{{\mathbb Z}}
\def\R{{\mathbb R}}

\def\eps{{\varepsilon}}
\def\e{\textbf{e}}

\def\mand{\qquad \mbox{and} \qquad}

\def\\{\cr}
\def\({\left(}
\def\){\right)}
\def\[{\left[}
\def\]{\right]}
\def\<{\langle}
\def\>{\rangle}
\def\fl#1{\left\lfloor#1\right\rfloor}
\def\rf#1{\left\lceil#1\right\rceil}
\def\le{\leqslant}
\def\ge{\geqslant}

\begin{document}

\title{Character sums with Beatty sequences\break on Burgess-type intervals}

\author{
{\sc William D.~Banks} \\
{Department of Mathematics} \\
{University of Missouri} \\
{Columbia, MO 65211 USA} \\
{\tt bbanks@math.missouri.edu} \\
\and
{\sc Igor E.~Shparlinski} \\
{Department of Computing}\\
{Macquarie University} \\
{Sydney, NSW 2109, Australia} \\
{\tt igor@ics.mq.edu.au}}

\maketitle

\begin{abstract}
We estimate multiplicative character sums taken on the values of a
non-homogeneous Beatty sequence $\{\fl{\alpha n + \beta}~:~n
=1,2,\ldots\,\}$, where $\alpha,\beta\in\R$, and $\alpha$ is
irrational. Our bounds are nontrivial over the same short intervals
for which the classical character sum estimates of Burgess have been
established.
\end{abstract}


\section{Introduction}

For fixed $\alpha,\beta\in\R$, the corresponding
\emph{non-homogeneous Beatty sequence} is the sequence of integers
defined by
$$
\cB_{\alpha,\beta}=\(\fl{\alpha n+\beta}\)_{n=1}^\infty,
$$
where $\fl{x}$ denotes the greatest integer $\le x$ for every
$x\in\R$. Beatty sequences arise in a variety of apparently
unrelated mathematical settings, and because of their versatility,
the arithmetic properties of these sequences have been extensively
explored in the literature; see, for example,
\cite{Abe,Beg,Kom1,Kom2,OB,Tijd} and the references contained
therein.

In this paper, we study character sums of the form
$$
S_k(\alpha,\beta,\chi;N)=\sum_{n\le N}\chi(\fl{\alpha n+\beta}),
$$
where $\alpha$ is \emph{irrational}, and $\chi$ is a
\emph{non-principal character modulo~$k$}. In the special case that
$k=p$ is a prime number, the sums $S_p(\alpha,\beta,\chi;N)$ have
been previously studied and estimated nontrivially for $N\ge
p^{1/3+\eps}$, where $\varepsilon>0$; see~\cite{BaSh1,BaSh2}.

Here, we show that the approach of~\cite{Abe} (see
also~\cite{Beg,LuZh}), combined with a bound on sums of the form
\begin{equation}
\label{eq:Sums U} U_k(t,\chi;M_0,M)=\sum_{M_0< m\le M}\chi(m) \,\e(t
m)\qquad(t\in\R),
\end{equation}
where $\e(x) = \exp(2 \pi i x)$ for all $x\in\R$, yields a
nontrivial bound on the sums $S_k(\alpha,\beta,\chi;N)$ for all
sufficiently large $N$ (see Theorem~\ref{thm:char} below for a
precise statement). In particular, in the case that $k=p$ is prime,
we obtain a nontrivial bound for all $N\ge p^{1/4+\eps}$, which
extends the results found in~\cite{BaSh1,BaSh2}.

It has recently been shown in~\cite{BaGaSh} that for a prime $p$ the
least positive quadratic non-residue modulo~$p$ among the terms of a
Beatty sequence is of size at most $p^{1/(4e^{1/2})+ o(1)}$, a
result which is complementary to ours.  However, the underlying
approach of~\cite{BaGaSh} is very different and cannot be used to
bound the sums $S_k(\alpha,\beta,\chi;N)$.

We remark that one can obtain similar results to ours by using
bounds for double character sums, such as those given
in~\cite{FrIw}. The approach of this paper, however, which dates
back to~\cite{Abe}, seems to be more general and can be used to
estimate similar sums with many other arithmetic functions $f(m)$
provided that appropriate upper bounds for the sums
$$
V(t,f;M_0,M)=\sum_{M_0< m\le M}f(m) \,\e(t m)\qquad(t\in\R)
$$
are known. Such estimates have been obtained for the characteristic
functions of primes and of smooth numbers (see~\cite{Dav}
and~\cite{FoTe}, respectively), as well as for many other functions.
Thus, in principle one can obtain asymptotic formulas for the number
of primes or smooth numbers in a segment of a Beatty sequence (in
the case of smooth numbers, this has been done in~\cite{BaSh3} by a
different method).

\bigskip

\noindent{\bf Acknowledgements.} The authors would like to thank
Moubariz Garaev for several fruitful discussions. This work began
during a pleasant visit by W.~B.\ to Macquarie University; the
support and hospitality of this institution are gratefully
acknowledged. During the preparation of this paper, I.~S.\ was
supported in part by ARC grant DP0556431.

\section{Notation}

Throughout the paper, the implied constants in the symbols $O$ and
$\ll$ may depend on $\alpha$ and $\eps$ but are absolute otherwise.
We recall that the notations $U = O(V)$ and $U\ll V$ are equivalent
to the assertion that the inequality $|U|\le c\,V$ holds for some
constant $c>0$.

We also use the symbol $o(1)$ to denote a function that tends to $0$
and depends only on $\alpha$ and $\eps$. It is important to note
that our bounds are uniform with respect all of the involved
parameters other than $\alpha$ and $\eps$; in particular, our bounds
are uniform with respect to $\beta$. In particular, the latter means
that our bounds also apply to the shifted sums of the form
$$
\sum_{M+1 \le n\le M+ N}\chi(\fl{\alpha n+\beta}) = \sum_{n\le
N}\chi(\fl{\alpha n+ \alpha M + \beta}),
$$
and these bounds are uniform for all integers $M$.

In what follows, the letters $m$ and $n$ always denote non-negative
integers unless indicated otherwise.

We use $\fl{x}$ and $\{x\}$ to denote the greatest integer $\le x$
and the fractional part of $x$, respectively.

Finally, recall that the \emph{discrepancy} $D(M)$ of a sequence of
(not necessarily distinct) real numbers $a_1,\ldots,a_M\in[0,1)$ is
defined by
\begin{equation}
\label{eq:descr_defn}
D(M)=\sup_{\cI\subseteq[0,1)}\left|\frac{V(\cI,M)}{M}-|\cI|\,\right|,
\end{equation}
where the supremum is taken all subintervals $\cI=(c,d)$ of the
interval $[0, 1)$, $V(\cI,M)$ is the number of positive integers
$m\le M$ such that $a_m\in\cI$, and $|\cI|=d-c$ is the length of
$\cI$.

\section{Preliminaries}

It is well known  that for every irrational number $\alpha$, the
sequence of fractional parts
$\{\alpha\},\{2\alpha\},\{3\alpha\},\ldots\,$, is \emph{uniformly
distributed modulo~$1$} (for instance, see \cite[Example~2.1,
Chapter~1]{KuNi}). More precisely, let $D_{\alpha,\beta}(M)$ denote
the discrepancy of the sequence $(a_m)_{m=1}^M$, where
$$
a_m=\{\alpha m + \beta\}\qquad(m=1,2,\ldots,M).
$$
Then, we have:

\begin{lemma}
\label{le:discr irrat}  Let $\alpha$ be a fixed irrational number.
Then, for all $\beta\in\R$ we have
$$
D_{\alpha,\beta}(M)\le 2D_{\alpha,0}(M)=o(1)\qquad(M\to\infty),
$$
where the function implied by $o(1)$ depends only on $\alpha$.
\end{lemma}

When more information about $\alpha$ is available, the bound of
Lemma~\ref{le:discr irrat} can be made more explicit. For this, we
need to recall some familiar notions from the theory of
\emph{Diophantine approximations}.

For an irrational number $\alpha$, we define its \emph{type} $\tau$
by the relation
$$
\tau=\sup\Bigl\{\vartheta\in\R~:~\liminf_{q\to\infty,~q\in\Z^+}
q^\vartheta\,\|\alpha q\|=0\Bigr\}.
$$
Using~\emph{Dirichlet's approximation theorem}, it is easy to see
that $\tau\ge 1$ for every irrational number $\alpha$. The
celebrated theorems of Khinchin~\cite{Khin} and of Roth~\cite{Roth}
assert that $\tau=1$ for almost all real numbers $\alpha$ (with
respect to Lebesgue measure) and all algebraic irrational numbers
$\alpha$, respectively; see also~\cite{Bug,Schm}.

The following result is taken from~\cite[Theorem~3.2,
Chapter~2]{KuNi}:

\begin{lemma}
\label{le:discr with type}  Let $\alpha$ be a fixed irrational
number of type $\tau<\infty$.  Then, for all $\beta\in\R$ we have
$$
D_{\alpha,\beta}(M)\le M^{-1/\tau+o(1)}\qquad(M\to\infty),
$$
where the function implied by $o(1)$ depends only on $\alpha$.
\end{lemma}

Next, we record the following property of type:

\begin{lemma}
\label{le:type stability} If $\alpha$ is an irrational number of
type $\tau<\infty$ then so are $\alpha^{-1}$ and $a \alpha$
   for any integer $a \ge 1$.
\end{lemma}

Finally, we need the following elementary result, which describes
the set of values taken by the Beatty sequence $\cB_{\alpha,\beta}$
in the case that $\alpha>1$:
\begin{lemma}
\label{lem:Beatty values} Let $\alpha > 1$.  An  integer $m$ has the
form $m= \fl{\alpha n + \beta}$ for some integer $n$ if and only if
$$
0<\{\alpha^{-1}(m-\beta+1)\}\le\alpha^{-1}.
$$
The value of $n$ is determined uniquely by $m$.
\end{lemma}

\begin{proof}
  It is easy to see that an integer
$m$ has the form $m=\fl{\alpha n+\beta}$ for some integer $n$ if and
only if the inequalities
$$
\frac{m-\beta}{\alpha} \le n < \frac{m-\beta+1}{\alpha}
$$
hold, and since $\alpha>1$ the value of $n$ is determined uniquely.
\end{proof}

\section{Character Sums}
\label{sec:char sums}

For every real number $\eps>0$ and integer $k\ge 1$, we put
\begin{equation}
\label{eq:Bdefine}
B_\eps(k)=\left\{%
\begin{array}{ll}
      k^{1/4+\eps} & \quad\hbox{if $k$ is prime;} \\
      k^{1/3+\eps} & \quad\hbox{if $k$ is a prime power;} \\
      k^{3/8+\eps} & \quad\hbox{otherwise.}
\end{array}%
\right.
\end{equation}

\begin{theorem}
\label{thm:char} Let $\alpha>0$ be a fixed irrational number, and
let $\eps>0$ be fixed. Then, uniformly for all $\beta\in\R$, all
non-principal multiplicative characters $\chi$ modulo $k$, and all
integers $N\ge B_\eps(k)$, we have
$$
S_k(\alpha,\beta,\chi;N)=o(N)\qquad(k\to\infty),
$$
where the function implied by $o(N)$ depends only on $\alpha$ and
$\eps$.
\end{theorem}

\begin{proof}
We can assume that $\eps<1/10$, and this implies that $B_\eps(k)\le
k^{2/5}$ in all cases. Observe that it suffices to prove the result
in the case that $B_\eps(k)\le N\le k^{1/2}$. Indeed, assuming this
has been done, for any $N>k^{1/2}$ we put $N_0=\fl{k^{9/20}}$ and
$t=\fl{N/N_0}$; then, since $B_\eps(k)\le N_0\le k^{1/2}$ we have
\begin{equation*}
\begin{split}
S_k(\alpha,\beta,\chi;N)&=\sum_{j=0}^{t-1}\sum_{n\le
N_0}\chi(\fl{\alpha(n+jN_0)+\beta}) +\sum_{tN_0<n\le
N}\chi(\fl{\alpha n+\beta})\\
&=\sum_{j=0}^{t-1}S_k(\alpha,\beta+\alpha jN_0,\chi;N_0)+O(N_0)\\
&=o(tN_0)+O\bigl(Nk^{-1/20}\bigr)=o(N)\qquad(k\to\infty)
\end{split}
\end{equation*}
using the fact that our bounds are uniform with respect to $\beta$.

We first treat the case that $\alpha>1$. Put $\gamma = \alpha^{-1}$,
$\delta = \alpha^{-1}(1-\beta)$, $M_0=\fl{\alpha+\beta-1}$, and
$M=\fl{\alpha N+\beta}$. {From} Lemma~\ref{lem:Beatty values} we see
that
\begin{equation}
\label{eq:Char Fun} S_k(\alpha,\beta,\chi;N)=\sum_{\substack{M_0<
m\le M\\ 0<\{\gamma m + \delta\}\le \gamma}}\chi(m) = \sum_{M_0<
m\le M}\chi(m) \,\psi(\gamma m + \delta),
\end{equation}
where $\psi(x)$ is the periodic function with period one for which
$$
\psi(x) = \left\{  \begin{array}{ll}
1& \quad \hbox{if $0< x\le \gamma$}; \\
0& \quad \mbox{if $\gamma<x\le 1$}.
\end{array} \right.
$$

By a classical result of Vinogradov (see~\cite[Chapter~2,
Lemma~2]{Vin}) it is known that for any $\Delta$ such that
$$
0 < \Delta < \frac{1}{8} \mand \Delta\le
\frac{1}{2}\,\min\{\gamma,1-\gamma\},
$$
there is a real-valued function $\psi_\Delta(x)$ with the following
properties:
\begin{itemize}
\item $\psi_\Delta(x)$ is periodic with period one;

\item $0 \le \psi_\Delta(x) \le 1$ for all $x\in\R$;

\item $\psi_\Delta(x) = \psi(x)$ if $\Delta\le x\le
\gamma-\Delta$ or $\gamma+\Delta\le x\le 1-\Delta$;

\item $\psi_\Delta(x)$ can be represented as a Fourier series
$$
\psi_\Delta(x) = \gamma + \sum_{j=1}^\infty \(\,g_j\,\e(jx) + h_j
\,\e(-jx)\),
$$
where the coefficients $g_j,h_j$ satisfy the uniform bound
$$
\max\{|g_j|, |h_j|\}\ll\min\{j^{-1},j^{-2}\Delta^{-1}\}
   \qquad(j\ge 1).
$$
\end{itemize}
Therefore, from~\eqref{eq:Char Fun} we derive that
\begin{equation}
\label{eq:Char Fun Approx} S_k(\alpha,\beta,\chi;N) = \sum_{M_0<
m\le M}\chi(m) \psi_\Delta(\gamma m + \delta)+O(V(\cI,M_0,M)),
\end{equation}
where $V(\cI,M_0,M)$ denotes the number of integers $M_0< m\le M$
such that
$$
\{\gamma m+\delta\}\in\cI=
[0,\Delta)\cup(\gamma-\Delta,\gamma+\Delta) \cup(1-\Delta,1).
$$
Since $|\cI|\ll\Delta$, it follows from Lemma~\ref{le:discr irrat}
and the definition~\eqref{eq:descr_defn} that
\begin{equation}
\label{eq:bound V(I,M)} V(\cI,M_0,M)\ll\Delta N+o(N),
\end{equation}
where the implied function $o(N)$ depends only on $\alpha$.

To estimate the sum in~\eqref{eq:Char Fun Approx}, we insert the
Fourier expansion for $\psi_\Delta(\gamma m + \delta)$ and change
the order of summation, obtaining
\begin{equation*}
\begin{split}
\sum_{M_0<m\le M}&\chi(m)\psi_\Delta(\gamma m + \delta)
=\gamma\,U_k(0,\chi;M_0,M)\\
&+\sum_{j=1}^\infty g_j\,\e(\delta j)\,U_k(\gamma
j,\chi;M_0,M)+\sum_{j=1}^\infty h_j\,\e(-\delta j)\,U_k(-\gamma
j,\chi;M_0,M),
\end{split}
\end{equation*}
where the sums $U_k(t,\chi;M_0,M)$ are defined by~\eqref{eq:Sums U}.

Since $M-M_0\ll N$, using the well known results of
Burgess~\cite{Burg1,Burg2,Burg3} on bounds for partial Gauss sums,
it follows that for any fixed $\varepsilon>0$ there exists $\eta>0$
such that
\begin{equation}
\label{eq:Burg} U_k(a/k,\chi;M_0,M) \ll N^{1-\eta}
\end{equation}
holds uniformly for all $N\ge B_\eps(k)$ and all integers $a$;
clearly, we can assume that $\eta\le 1/10$.

Put $r=\fl{\gamma k}$. Then, for any integer $n$, we have
$$
\e(\gamma n) - \e(rn/k)  \ll |\gamma n  -  rn/k| \le |n|k^{-1},
$$
which implies that
$$
U_k(\gamma j,\chi;M_0,M)  = U_k(rj/k,\chi;M_0,M) + O(N^2k^{-1}|j|).
$$
Using~\eqref{eq:Burg} in the case that $|j|\le kN^{-1-\eta}$ we
derive that
$$
U_k(\gamma j,\chi;M_0,M)\ll N^{1-\eta},
$$
and for $|j|>kN^{-1-\eta}$ we use the trivial bound
$$
\bigl|U_k(\gamma j,\chi;M_0,M)\bigr|\ll N.
$$
Consequently,
\begin{equation*}
\begin{split}
\sum_{M_0<m\le M}\chi(m)&\psi_\Delta(\gamma m + \delta)\\
&\ll N^{1-\eta}\sum_{j \le kN^{-1-\eta}}(|g_j| + |h_j|)+
N\sum_{j> kN^{-1-\eta}}(|g_j| + |h_j|) \\
&\ll N^{1-\eta} \sum_{j \le kN^{-1-\eta}}j^{-1}+
N \Delta^{-1} \sum_{j>  kN^{-1-\eta}}j^{-2} \\
&\ll N^{1-\eta} \log k + N^{2+\eta} \Delta^{-1}  k^{-1}.
\end{split}
\end{equation*}
Since $N^2\le k\le N^4$, we see that
\begin{equation}
\label{eq:boundthesum} \sum_{M_0<m\le M}\chi(m)\psi_\Delta(\gamma m
-\delta)\ll N^{1-\eta} \log N + N^{\eta} \Delta^{-1}.
\end{equation}

Inserting the bounds~\eqref{eq:bound V(I,M)}
and~\eqref{eq:boundthesum} into~\eqref{eq:Char Fun Approx}, choosing
$\Delta=N^{(\eta-1)/2}$, and taking into account that $0<\eta\le
1/10$, we complete the proof in the case that $\alpha>1$.

If $\alpha<1$, put $a=\rf{\alpha^{-1}}$ and write
\begin{equation*}
\begin{split}
S_k(\alpha,\beta,\chi;N)&=\sum_{n\le N}\chi(\fl{\alpha n+\beta})\\
&=\sum_{j =0}^{a-1} \sum_{m \le (N-j)/a}
\chi(\fl{\alpha a m + \alpha j +\beta})\\
&=\sum_{j =0}^{a-1} S_k(\alpha a,\alpha j+\beta,\chi;(N-j)/a).
\end{split}
\end{equation*}
Applying the preceding argument with the irrational number $\alpha
a>1$, we conclude the proof.
\end{proof}

For an irrational number $\alpha$ of type $\tau<\infty$, we proceed
as in the proof of Theorem~\ref{thm:char}, using Lemma~\ref{le:discr
with type} instead of Lemma~\ref{le:discr irrat}, and also applying
Lemma~\ref{le:type stability}; this yields the following statement:

\begin{theorem}
\label{thm:char with type} Let $\alpha>0$ be a fixed irrational
number of type $\tau<\infty$. For every fixed $\eps>0$ there exists
$\rho>0$, which depends only on $\eps$ and $\tau$, such that for all
$\beta\in\R$, all non-principal multiplicative characters $\chi$
modulo $k$, and all integers $N\ge B_\eps(k)$, we have
$$
S_k(\alpha,\beta,\chi;N) \ll N k^{-\rho}.
$$
\end{theorem}

\end{document}